\documentclass[12pt]{article}

\usepackage{bm}
\usepackage{amsmath}
\usepackage{graphicx,times}
\usepackage[colorlinks=true, pdfstartview=FitV, linkcolor=black, citecolor= blue, urlcolor= black]{hyperref}

\renewcommand{\Vec}[1]{\mbox{\boldmath$#1$}}

\begin{document}

\title{Short-axis-mode rotation of a free rigid body by perturbation series}

\author{Martin Lara\thanks{Columnas de H\'ercules 1, ES-11100 San Fernando, Spain, \tt mlara0@gmail.com} 
}

\maketitle{}

\begin{abstract}
A simple rearrangement of the torque free motion Hamiltonian shapes it as a perturbation problem for bodies rotating close to the principal axis of maximum inertia, independently of their triaxiality. The complete reduction of the main part of this Hamiltonian via the Hamilton-Jacobi equation provides the action-angle variables that ease the construction of a perturbation solution by Lie transforms. The lowest orders of the transformation equations of the perturbation solution are checked to agree with Kinoshita's corresponding expansions for the exact solution of the free rigid body problem. For approximately axisymmetric bodies rotating close to the principal axis of maximum inertia, the common case of major solar system bodies, the new approach is advantageous over classical expansions based on a small triaxiality parameter.
\end{abstract}

\section{Introduction}

The rotational dynamics of solar system bodies is fundamentally driven by its torque free motion, although the tiny influence of external torques like the gravity-gradient or the YORP effect may induce long-period or secular effects. Nevertheless, the fact that these torques are generally small makes the computation of approximate analytical solutions possible using perturbation methods (see, for instance, \cite{Kinoshita1977,GetinoFerrandiz1991,CicaloScheeres2010}).
\par

The torque free motion admits analytical, closed form solution, which is commonly taken as the zero order in the computation of approximate solutions of perturbed problems. In a Hamiltonian framework, this solution is obtained by finding the transformation from Andoyer to action-angle variables \cite{Sadov1970,Kinoshita1972}, mainly because action-angle variables make achievable the essential checking of KAM conditions in the completely reduced Hamiltonian \cite{Kozlov1974}, and prevent the appearance of mixed secular-periodic terms in the solution of perturbed problems \cite{Kinoshita1972,LaraFerrer2013}. Then, the disturbing function is expressed as a trigonometric series in action-angle variables so that the perturbed Hamiltonian can be simplified by the usual averaging of fast-evolving angles \cite{Arnold1989}.
\par

The perturbation approach is quite feasible in the case of axisymmetric bodies, an instance in which Andoyer variables are action-angle variables by themselves. In that case, even the standard successive approximations method is suitable for finding the direct integration of the equations of motion by simple quadratures (see, for instance, \cite{LaraFukushimaFerrer2010} and references therein). In the triaxial case, however, because of the elliptic integrals and functions on which the solution of the unperturbed problem rests, expanding the disturbing function as a trigonometric series may require a nontrivial preprocessing in which elliptic integrals and functions are conveniently expressed in terms of Jacobi theta functions \cite{Sadov1970,Barkin1998}. Notwithstanding, some exceptions are found where the solution of the perturbed problem can be computed directly in elliptic functions, at least up to the first order, without need of resorting to series expansions \cite{Chernousko1963,HitzlBreakwell1971,LaraFerrer2013}.
\par

Alternatively, it is customary to expand the solution of the torque free motion in terms of certain small quantity. Thus, either in the (physical) case of rigid bodies with small triaxiality, or in the (dynamical) case of rotational motion close to a principal axis of inertia, early truncations of the series expansion of the unperturbed solution may result very efficient \cite{Kinoshita1992,SouchayFolgueiraBouquillon2003}. Specifically, the former applies to the rotation of the Earth and other major bodies of the solar system, whose figure departs only slightly from the axisymmetrical case \cite{Kinoshita1972,Kinoshita1977,GetinoEscapaMiguel2010,CottereauSouchayAljbaae2010}, whereas the latter may be the alternative for studying the rotation of some asteroids \cite{SouchayKinoshitaNakaiRoux2003,SouchayBouquillon2005}.
\par

It worths noting, however, that when the triaxiality of the rigid body is small, the part of the torque free rotation Hamiltonian related to triaxiality can be taken as a perturbation of the axisymmetric case, and, consequently, can be included among the terms of the disturbing function (see \cite{Zanardi1986,FerrandizSansaturio1989}, for instance). When applying this strategy to the integrable problem of torque free rotation, the truncated expansions of Kinoshita's analytical solution\cite{Kinoshita1972}  are recovered term for term \cite{FerrerLara2010}. This alternative approach, taking an integrable problem as a perturbation of a simpler integrable problem, may be applied to further rearrange the torque free motion Hamiltonian like a perturbed spherical rotor. Such an ordering makes a perturbation approach by Lie transforms very efficient in the case of almost spherical, almost axisymmetrical rigid bodies \cite{LaraFukushimaFerrer2011}. 
\par

The perturbation approach is not limited to the case of rigid bodies with small triaxiality and may be applied also to bodies with any triaxiality, either large or small, when the rotation happens close to the principal axis of maximum inertia, that is: rotating in the so-called short-axis-mode, or SAM (see Appendix C of \cite{Fukushima2008} for a detailed classification of rotational modes). The demonstration of the feasibility of such perturbation approach is the goal of the present research.
\par

The first step in the procedure is to reorganize the torque free Hamiltonian in such a way that the sine of the inclination angle between the body's equatorial plane and the plane perpendicular to the angular momentum vector is shown to scale the Hamiltonian terms related to the non-sphericity of the rigid body. Next, a trivial manipulation changes this general scaling into a perturbation scaling in which the magnitude of the small parameter is proportional to the sine of \emph{half} that inclination angle. Remarkably, the main part of the new perturbation arrangement, which is called here the main problem of SAM rotation, is formally the same as the torque free rotation Hamiltonian except for the halving of the inclination angle. But this subtle difference is enough to release the solution from its dependence on elliptic functions and integrals.
\par

It follows the standard Hamilton-Jacobi reduction \cite{Arnold1989} to find the action-angle variables of the main problem of SAM rotation, which is selected as the unperturbed problem. Then, after formulating the perturbation part of the free rigid body Hamiltonian in the new action-angle variables, the straightforward computation of a perturbation solution to the torque free motion by Lie transforms \cite{Hori1966,Deprit1969} is carried out up to any desired order. Comparison of the lowest orders of the perturbation series solution with Kinoshita's expansions of the closed form analytical solution of the free rigid body \cite{Kinoshita1992} shows complete agreement between both independent solutions.
\par

Finally, the known equivalence between short- and long-axis-mode (LAM) when interchanging the moments of maximum and minimum inertia, trivially shows that the perturbation solution is useful also in the case of LAM rotation with minor modifications in the definition of the inertia parameters.

\section{Hamiltonian arrangement}

Because of the preservation of the angular momentum vector $\Vec{M}$, the plane perpendicular to it (the so-called \emph{invariable} plane) plays an important role in finding the analytical solution of the torque free motion. Thus, classical solutions using Euler angles take this plane as the inertial plane to which the attitude of the rigid body is referred \cite{Whittaker1917,Golubev1960}. Referring the body frame to a different inertial plane is achieved by trivial rotations of fixed angles, say a new set of Euler angles.

This description of the attitude dynamics by means of two sets of Euler angles makes it natural the introduction of the angles $\lambda$, between the $x$ axis of the inertial plane and the ascending node of the invariable plane on the inertial plane, $\mu$, between this node and the ascending node of the equatorial plane of the body on the invariable plane, and $\nu$, between the latter node and the $x$ axis of the body. Their conjugate momenta are the modulus of the angular momentum vector $M=\|\Vec{M}\|$, conjugate to $\mu$, and the projections of $\Vec{M}$ over the $z$ axis of the inertial and body frames, $\Lambda$ and $N$, conjugate to $\lambda$ and $\nu$ respectively. The canonical set $(\lambda,\mu,\nu,\Lambda,M,N)$ is what is customarily called Andoyer variables \cite{Andoyer1923}, and so it is throughout this paper.
\par

The convenience of Andoyer variables results evident when using Hamiltonian formulation, because they disclose the integrable character of the torque free motion. Due to the absence of external forces the Hamiltonian of the torque free motion coincides with the kinetic energy of rotation about the center of mass $T=\frac{1}{2}\Vec\omega\,\mathbf{I}\,\Vec\omega$, where $\Vec\omega$ is the angular velocity and $\mathbf{I}$ is the inertia tensor. Therefore, taking into account that $\Vec{M}=\mathbf{I}\,\Vec\omega$ and choosing the body frame as defined by the principal axes of inertia, it is easily found that
\begin{equation}\label{Andoyer}
\mathcal{H}_0=\left(\frac{\sin^2\nu}{A}+\frac{\cos^2\nu}{B}\right)\frac{M^2-N^2}{2}+\frac{N^2}{2C},
\end{equation}
where $A\le{B}\le{C}$ are the body's principal moments of inertia.
\par

Equation (\ref{Andoyer}) shows that $\lambda$, $\Lambda$ and $M$ are integrals of the problem due to the ignorable character of their respective conjugate variables. In consequence, the Hamilton equations of Eq.~(\ref{Andoyer}) decouple the reduced system in $(\nu,N)$, of one degree of freedom, from the time evolution of $\mu$. The latter is solved by quadrature after integrating the reduced system. Therefore, the reduced phase space $(\nu,N)$ can be represented by simple contour plots of the Hamiltonian, providing an alternative description to the Poinsot ellipsoid representation \cite{DepritAJP1967}.
\par

Concerned with the rotation of the major bodies of the solar system, for which $C/A\approx C/B\approx1$, Andoyer\cite{Andoyer1923} introduced the dimensionless parameters $\alpha\ge0$ and $0\le\beta\le1$ given by
\begin{equation} \label{alfabeta}
\alpha\,(1+\beta)=\frac{C}{A}-1,
\qquad
\alpha\,(1-\beta)=\frac{C}{B}-1,
\end{equation}
and reordered Eq.~(\ref{Andoyer}) as a perturbation problem
\begin{equation}\label{hacheA}
\mathcal{H}_0=\frac{M^2}{2C}\left[1+\alpha\left(1-\frac{N^2}{M^2}\right)-\alpha\,\beta\left(1-\frac{N^2}{M^2}\right)\cos2\nu\right].
\end{equation}
\par

Both Andoyer's parameters have a clear physical meaning as is apparent in Eq.~(\ref{alfabeta}). Thus, $\beta=0$ implies $A=B$, the oblate case, whereas $\beta=1$ results from $C=B$, the prolate one, both extremes corresponding to axisymmetry. Therefore, $\beta$ shows how much the body's figure departs from the axisymmetrical case, and hence is usual called the triaxiality coefficient. On the other side, $\alpha=0$ corresponds to a spherical rotor $A=B=C$ and, therefore, this coefficient indicates how much the body departs from the spherical figure.\footnote{Andoyer's triaxiality coefficient $\beta$ is noted $e$ by Kinoshita \cite{Kinoshita1972} who, in addition, instead of $\alpha$ uses the inertia parameter $D=-C/\alpha$ first proposed by Hitzl and Breakwell \cite{HitzlBreakwell1971}.}
\par

Andoyer's arrangement of the free rigid body Hamiltonian is quite convenient when studying the rotation of the major bodies of the solar system under external torques, a case in which both $\alpha$ and $\beta$ are small quantities. Besides, this is the natural ordering in a perturbed spherical rotor approach, in which both the non sphericity and the triaxiality are taken as perturbations, and that proved very efficient when approaching the perturbed rotation of almost spherical bodies by Lie transforms \cite{FerrerLara2010,LaraFukushimaFerrer2011}.
\par

However, neither $\alpha$ nor $\beta$ have to be small, in general, and hence the general solution of the torque free motion must be taken as the zero order in the perturbation approach, a fact that complicates the development of the disturbing function because of the elliptic integrals and functions on which the solution of the free rigid body rotation unavoidably depends. Nevertheless, using the general solution to the torque free motion can be avoided in other specific cases, where perturbation strategies can be set up.
\par

That is the case of rotations close to the principal axis of maximum inertia, a case in which $N\approx{M}$. Then, the inclination angle $J$ between the equatorial plane of the rigid body and the invariable plane is small, and so it is $(1-N^2/M^2)=\sin^2\!J$ in Eq.~(\ref{hacheA}). However, because $M$ is an integral of the torque free motion,
the sine of the inclination just scales the reduced problem, and the form in which Eq.~(\ref{hacheA}) is written, namely, 
\begin{equation}\label{hacheAsin}
\mathcal{H}_0=\frac{M^2}{2C}\left[1+\alpha\sin^2J\left(1-\beta\cos2\nu\right)\right],
\end{equation}
is not amenable for a perturbation approach. Nevertheless, based on the standard relation
\[
\sin^2J=2\sin^2\frac{J}2\left(2-2\sin^2\frac{J}2\right),
\] 
trivial manipulations allow to cast Eq.~(\ref{hacheA}) in the form of a perturbation problem
\begin{equation}\label{hacheP}
\mathcal{H}_0=\mathcal{A}+\varepsilon\,\mathcal{P},
\end{equation}
where
\begin{equation}\label{main}
\mathcal{A}=\frac{M^2}{2C}\left[1+2\alpha\,\left(1-\frac{N}{M}\right)\left(1-\beta\cos2\nu\right)\right],
\end{equation}
is the integrable part, $\varepsilon$ is a formal small parameter, and
\begin{equation}\label{disturbing}
\mathcal{P}=-\frac{M^2}{2C}\,\alpha\left(1-\frac{N}{M}\right)^2\left(1-\beta\cos2\nu\right),
\end{equation}
is the perturbation. In fact, we should have written $|N|$ in Eqs.~(\ref{main}) and (\ref{disturbing}) instead of $N$. However, both cases are symmetric and hence we limit our discussion to the case $N>0$.
\par

Since Eq.~(\ref{hacheP}) is a perturbation problem independently of the value of $\beta$, this perturbation Hamiltonian can be used to represent the torque free rotation of major bodies of the solar system as well as asteroids rotating close to the minimum energy state. Because of that, in what follows Eq.~(\ref{main}) is called the main problem of SAM rotation, or main problem in short.

\section{The main problem of SAM rotation}

When Eq.~(\ref{main}) is rewritten as a function of the inclination angle $J$, it takes the form
\begin{equation}\label{mainJ}
\mathcal{A}=\frac{M^2}{2C}\left[1+\alpha'\sin^2(J/2)\left(1-\beta\cos2\nu\right)\right],
\qquad
\alpha'=4\alpha,
\end{equation}
which is formally equal to Eq.~(\ref{hacheAsin}), except for each Hamiltonian deals with a different inclination angle: $\gamma=J$ in the case of the torque free motion and $\gamma=J/2$ for the main problem of SAM rotation. However, the subtle difference introduced by the halving of $J$ in Eq.~(\ref{mainJ}) makes that the main problem can be integrated using trigonometric functions, as opposite to the elliptic functions and integrals required in the solution of Eq.~(\ref{hacheAsin}).
\par

\begin{figure}[htb]
\centering
\includegraphics[scale=1.3]{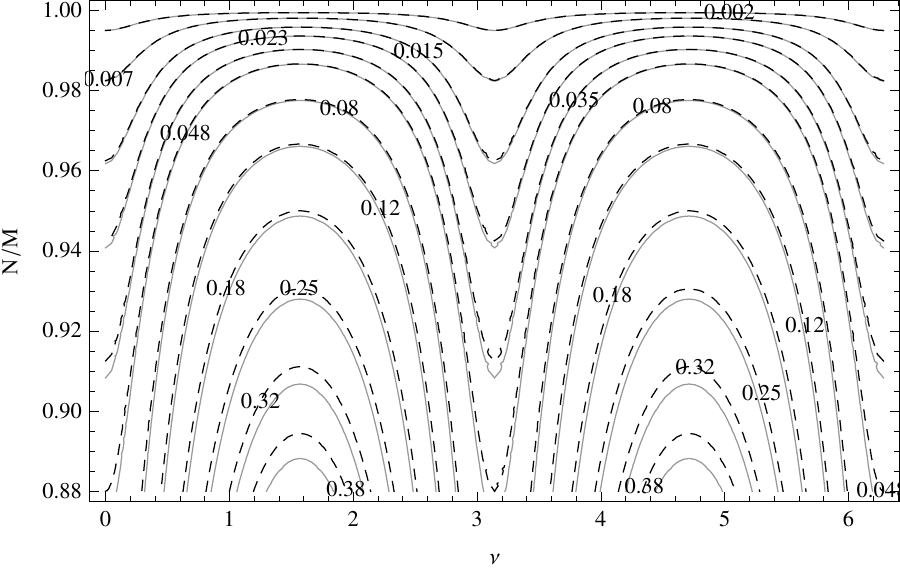} 
\caption{Reduced phase space of the free rigid body (full lines) and of the main problem of the asteroid rotation (dashed lines) for a triaxiality coefficient $\beta=0.8$ and $0.88\le{N}/M\le1$.}
\label{f:phases}
\end{figure}

Besides, the phase space of both different problems shows a nice agreement for small enough values of $J$, as illustrated in Fig.~\ref{f:phases}, where the flow in the reduced $(\nu,N)$ space is obtained from simple contour plots of the non-constant part of Hamiltonians (\ref{hacheAsin}) and (\ref{mainJ}) after scaling by $\alpha$. Namely, $\mathcal{H}_0^*=(2C\,\mathcal{H}_0/M^2-1)/\alpha$ and  $\mathcal{A}^*=(2C\,\mathcal{A}/M^2-1)/\alpha$, respectively. Specifically, Fig.~\ref{f:phases} presents contours $\mathcal{H}_0^*=\mathcal{A}^*=q$ for $q=0.002$, $0.007$, $0.015$, $0.023$, $0.035$, $0.048$, $0.08$, $0.12$, $0.18$, $0.25$, $0.32$, and $0.38$.

With the aim of using the main problem of SAM rotation as the zero order of a perturbation theory, the integration of the flow associated to Eq.~(\ref{mainJ}) is achieved by Hamiltonian reduction, using the Hamilton-Jacobi method, rather than approaching the direct integration of the Hamilton equations of $\mathcal{A}$.

\subsection{Hamilton-Jacobi reduction}

The essence of the Hamilton-Jacobi method is to find a canonical transformation such that the Hamiltonian flow in the new variables can be solved explicitly as a function of time. In consequence, the original flow is integrated by the simple expedient of introducing this time explicit solution into the transformation equations of the canonical transformation. One common approach is to find the transformation to action-angle variables because it makes cyclic all the coordinates in the Hamiltonian, which is, therefore, trivially integrated \cite{Arnold1989}.
\par

In the case of concern, a transformation $(\ell,g,L,G)\rightarrow(\mu,\nu,M,N)$ is found such that it converts Eq.~(\ref{main}) in a new Hamiltonian $\Phi\equiv\Phi(-,-,L,G)$ that depends only on the new momenta. Such canonical transformation is derived from a generating function in mixed variables $S\equiv{S}(\mu,\nu,L,G)$, depending on the old coordinates and the new momenta. The transformation equations are
\begin{equation}\label{tranequ}
(\ell,g,M,N)=\frac{\partial{S}}{\partial(L,G,\mu,\nu)}.
\end{equation}
\par

Because Eq.~(\ref{main}) does not depend on $\mu$, the generating function is chosen of the form
\[
S=\mu\,G+W(\nu,L,G).
\]
Therefore, Eq.~(\ref{tranequ}) is written in terms of the \emph{characteristic} function $W$ as the mixed transformation
\begin{equation}\label{treq}
M=G,\qquad 
N=\frac{\partial{W}}{\partial\nu},\qquad 
\ell=\frac{\partial{W}}{\partial{L}},\qquad 
g=\mu+\frac{\partial{W}}{\partial{G}}.
\end{equation}
\par

The Hamilton-Jacobi equation is then obtained by replacing in Eq.~(\ref{main}) the transformation equations of $M$ and $N$ given in Eq.~(\ref{treq}). Namely,
\begin{equation}\label{HJ}
\frac{G^2}{2C}\left[1+2\alpha\left(1-\frac{1}{G}\,\frac{\partial{W}}{\partial\nu}\right)\left(1-\beta\,\cos2\nu\right)\right]=\Phi(L,G).
\end{equation}
\par

Because the new Hamiltonian $\Phi$ in Eq.~(\ref{HJ}) does not depend on $\nu$, the selection of a particular form of $\Phi$ can be postponed (see \cite{FerrerLara2010b} and references therein). Anyway, it is found convenient to choose $\Phi$ of the form
\begin{equation}\label{standard}
\Phi=\frac{G^2}{2C}\left(1+2\alpha\,\frac{\Psi}{G}\right),
\end{equation}
where the function $\Psi\equiv\Psi(L,G)$ still remains to be determined.
\par

Equation (\ref{HJ}) is solved for $W$ to give
\begin{equation}\label{genw}
W=G\,\nu-\Psi\,\int\frac{1}{1-\beta\,\cos2\nu}\,\mathrm{d}\nu,
\end{equation}
where the quadrature is conveniently solved making the change of variable
\begin{equation}\label{nu2u}
\begin{array}{rcl}
\sin\nu &=& \displaystyle\frac{\sqrt{1-\beta}\,\sin{u}}{\sqrt{(1-\beta)\sin^2u+(1+\beta)\cos^2u}}, \\[4ex]
\cos\nu &=& \displaystyle\frac{\sqrt{1+\beta}\,\cos{u}}{\sqrt{(1-\beta)\sin^2u+(1+\beta)\cos^2u}},
\end{array}
\end{equation}
from which is derived the differential identity
\begin{equation}\label{differentials}
\mathrm{d}\nu=\frac{\sqrt{1-\beta^2}}{1+\beta\cos2u}\,\mathrm{d}u.
\end{equation}
\par

Note that Eq.~(\ref{nu2u}) is the standard change of variable used to solve the free rigid body rotation in closed form, c.f.~Eq.~(43) of \cite{HitzlBreakwell1971}. Note also that Eq.~(\ref{nu2u}) can be rewritten as
\[
\tan\nu=\sqrt{\frac{1-\beta}{1+\beta}}\tan{u},
\]
whose analogy with the relation between the true and eccentric anomalies in the Keplerian motion was emphasized in \cite{Kinoshita1972}.\par

Then, after standard manipulation, Eq.~(\ref{genw}) is solved to give
\[
W=G\,\nu-\frac{\Psi}{\sqrt{1-\beta^2}}\,u,
\]
which is replaced in the transformation equations (\ref{treq}) to give rise to a family of canonical transformations
\begin{eqnarray*}
M &=& G, \\
N &=& G-\frac{\Psi}{1-\beta\,\cos2\nu}=G-\Psi\,\frac{1+\beta\cos2u}{1-\beta^2}, \\
\ell &=& -\frac{\partial\Psi}{\partial{L}}\frac{u}{\sqrt{1-\beta^2}}, \\
g &=& \mu+\nu-\frac{\partial\Psi}{\partial{G}}\frac{u}{\sqrt{1-\beta^2}},
\end{eqnarray*}
from which specific transformations are derived for each particular selection of $\Psi$. Arbitrary choices of $\Psi$ provide, of course, different solutions for the integrable Hamiltonian (\ref{main}). However, if the solution is intended to be used in a perturbation approach, the suitability of each selection of $\Psi$, and consequently of the new Hamiltonian (\ref{standard}), should be carefully explored. 
\par

The usual requirement for the new variables to be action-angle variables, results in the natural choice
\begin{equation}\label{newham}
\Psi=\sqrt{1-\beta^2}\,L,
\end{equation}
and hence, the new, reduced Hamiltonian
\begin{equation}\label{reduced}
\Phi=\frac{G^2}{2C}\left(1+2\alpha\sqrt{1-\beta^2}\,\frac{L}{G}\right),
\end{equation}
which is quadratic in the momenta and whose Hessian determinant never vanishes.
\par

Therefore, the transformation from Andoyer variables to action-angle variables is
\begin{eqnarray} \label{scnu}
\sin\nu &=&-\frac{\sqrt{1-\beta}\,\sin\ell}{\sqrt{1+\beta\cos2\ell}},\qquad
\cos\nu=\frac{\sqrt{1+\beta}\,\cos\ell}{\sqrt{1+\beta\cos2\ell}}, \\[1ex] \label{mug}
\mu &=& g-\nu, \\
M &=& G, \\ \label{NL}
N &=& G-L\,\frac{1+\beta\cos2\ell}{\sqrt{1-\beta^2}},
\end{eqnarray}
whereas the inverse transformation is
\begin{eqnarray}\label{scl}
\sin\ell &=& -\frac{\sqrt{1+\beta}\,\sin\nu}{\sqrt{1-\beta\cos2\nu}},\qquad
\cos\ell=\frac{\sqrt{1-\beta}\,\cos\nu}{\sqrt{1-\beta\cos2\nu}}, \\
g &=& \mu+\nu, \\
G &=& M,\\
L &=& (M-N)\,\frac{1-\beta\,\cos2\nu}{\sqrt{1-\beta^2}}.
\end{eqnarray}
\par

The final check
\[
M\,\mathrm{d}\mu+N\,\mathrm{d}\nu=L\,\mathrm{d}\ell+G\,\mathrm{d}g,
\]
shows that, in addition, the canonical transformation from Andoyer variables to the action-angle variables of the main problem of SAM rotation is a Mathieu transformation.

\subsection{Disturbing function in new variables}

It is easy to check from Eq.~(\ref{scnu}) that 
\[
\cos2\nu=\frac{\beta+\cos2\ell}{1+\beta\cos2\ell},\qquad
\sin2\nu=\frac{\sqrt{1-\beta^2}\sin2\ell}{1+\beta\cos2\ell}.
\]
Then, the disturbing function Eq.~(\ref{disturbing}) is expressed in the new variables as
\begin{equation}\label{disturbing2}
\mathcal{P}=-\frac{L^2}{2C}\,\alpha\,(1+\beta\cos2\ell).
\end{equation}
\par

In consequence, the torque free Hamiltonian (\ref{hacheA}) is written in the new variables as
\begin{equation}\label{rigid2}
\mathcal{K}=\frac{G^2}{2C}\left[1+2\alpha\,\sqrt{1-\beta^2}\,\frac{L}{G}\left(1-\frac{L}{G}\,\frac{1+\beta\cos2\ell}{2\sqrt{1-\beta^2}}\right)\right],
\end{equation}
which is a perturbation problem in those cases in which
\begin{equation}\label{delta}
\delta=\frac{L}{G}\,\frac{1+\beta\cos2\ell}{\sqrt{1-\beta^2}}=2\sin^2(J/2)\ll1.
\end{equation}
For instance, $J=1,2.5,8$ and $25$ deg, results in $\delta\approx10^{-4},10^{-3},10^{-2}$ and $0.1$, respectively.
\par

The condition in Eq.~(\ref{delta}) is quite common for solar system bodies rotating in the short-axis-mode, even in the case of highly triaxial asteroids. Thus, for instance, using the accepted values in Table \ref{t:ssb}, Eq.~(\ref{delta}) provides $\delta=\mathcal{O}(10^{-13})$ for Mars, $\delta=\mathcal{O}(10^{-11})$ for the Earth, $\delta=\mathcal{O}(10^{-10})$ for the Moon, or $\delta=\mathcal{O}(10^{-8})$ for Eros.
\par

\begin{table} [htbp]
\label{t:ssb}
\centering
\begin{tabular}{@{}llllll@{}}
\multicolumn{1}{c}{Body} & \multicolumn{1}{c}{$A/C$} & \multicolumn{1}{c}{$B/C$} & \multicolumn{1}{c}{$\beta$} & \multicolumn{1}{c}{$J_0$} & \multicolumn{1}{c}{Ref.} \\
\hline
Mars & $0.9942917$ & $0.9949813$ & $0.0646316$ & $0.1''$ & \cite{SouchayFolgueiraBouquillon2003,DehantdeVironKaratekinVanHoolst2006} \\
Earth & $0.9967200$ & $0.9967222$ & $0.0003366$ & $1''$ & \cite{Fukushima2006} \\
Moon & $0.999368$ & $0.999601$ & $0.226105$ & $6.2''$ & \cite{Fukushima2006,NewhallWilliams1997} \\
Eros & $0.229427$ & $0.963754$ & $0.977853$ & $55''$ & \cite{SouchayKinoshitaNakaiRoux2003} \\
\hline
\end{tabular}
\caption{Inertia parameters and inclination angle for different solar system bodies.}
\end{table}

Therefore, for typical solar system bodies in SAM rotation, a perturbation solution in the new variables, with the small parameter $\epsilon\propto\delta$, would converge much faster than classical perturbation solutions based on a small triaxiality coefficient, with the small parameter $\epsilon\propto\beta$.

\section{Perturbation solution}

The formulation of the torque free Hamiltonian in the new variables, Eq.~ (\ref{rigid2}), makes it amenable to being solved with standard perturbation methods by Lie transforms \cite{Hori1966,Deprit1969}. More precisely, Deprit's algorithm is used here because it is specifically designed for the automatic computation of higher orders by machine \cite{Deprit1969}.
\par

Thus, a perturbation solution to Eq.~(\ref{rigid2}) is obtained by finding a canonical transformation of the Lie type
\begin{equation}\label{tx}
(\ell,g,L,G)\longrightarrow(\ell',g',L',G')
\end{equation}
from old to new, prime, variables, such that, up to certain truncation order, the transformed Hamiltonian depends only on the new momenta. The procedure is standard these days, and hence only the results are presented here. The interested reader can find all the relevant details in the original references or in modern textbooks on celestial mechanics.
\par

In this way, it has been obtained the long-term, averaged Hamiltonian
\begin{equation}\label{averaged}
\mathcal{T}=\frac{G'^2}{2C}\left[1+2\alpha\,\frac{L'}{G'}\,\sqrt{1-\beta^2}-\alpha\frac{L'^2}{G'^2}\left(1+\beta^2\,\sum_{i\ge1}\,\delta'^i\,q_i\right)\right],
\end{equation}
where
\begin{equation}\label{deltaprime}
\delta'=\frac{L'}{G'}\,\frac{1}{\sqrt{1-\beta^2}},
\end{equation}
and the triaxiality polynomials $q_i$ are given in Table \ref{t:polyham} up to $i=10$.
\par

\begin{table} [htbp]
\label{t:polyham}
\centering
\begin{tabular}{rllr}
\hline\noalign{\smallskip}
$i$ & \multicolumn{1}{c}{$q_i$} & \multicolumn{1}{c}{$q_i$} & $i$  \\
\noalign{\smallskip}\hline\noalign{\smallskip}
$1$ & $\frac{1}{2}$ & $\frac{13}{8}+\frac{10803}{256}\beta^2+\frac{630357}{4096}\beta^4+\frac{232505}{2048}\beta^6+\frac{27937}{2048}\beta^8$ & $10$ \\[1ex]
$2$ & $\frac{5}{8}$ & $\frac{3}{2}+\frac{1785}{64}\beta^2+\frac{70179}{1024}\beta^4+\frac{237339}{8192}\beta^6+\frac{36597}{32768}\beta^8$ & $9$ \\[1ex]
$3$ & $\frac{3}{4}+\frac{9}{32}\beta^2$ & $\frac{11}{8}+\frac{9009}{512}\beta^2+\frac{55583}{2048}\beta^4+\frac{44825}{8192}\beta^6$ & $8$ \\[1ex]
$4$ & $\frac{7}{8}+\frac{35}{32}\beta^2$ & $\frac{5}{4}+\frac{2675}{256}\beta^2+\frac{9305}{1024}\beta^4+\frac{4765}{8192}\beta^6$ & $7$ \\[1ex]
$5$ & $1+\frac{177}{64}\beta^2+\frac{45}{128}\beta^4$ & $\frac{9}{8}+\frac{2925}{512}\beta^2+\frac{2385}{1024}\beta^4$ & $6$ \\
\noalign{\smallskip}\hline
\end{tabular}
\caption{Triaxiality polynomials of the secular Hamiltonian in Eq.~(\protect\ref{averaged})}
\end{table}

The transformation equations of the averaging are $G=G'$ and
\begin{eqnarray}\label{ggp}
g &=& g'-\frac{L'}{G'}\,\sum_{i\ge1}\,\delta'^i\sum_{m=1,k}(-\beta)^mg_{i,m}\sin2m\ell' \\ \label{llp}
\ell &=& \ell'+\sum_{i\ge1}\,\delta'^i\sum_{m=1,i}(-\beta)^m\ell_{i,m}\sin2m\ell' \\ \label{NNp}
L &=& L'+L'\sum_{i\ge1}\,\delta'^i\left(\beta^2\,L_{i,0}-\sum_{m=1,k}(-\beta)^mL_{i,m}\cos2m\ell'\right)
\end{eqnarray}
where $k=(i+1)/2$ is an integer division, and the triaxiality polynomials $g_{i,m}$, $\ell_{i,m}$ and $L_{i,m}$ are given in Tables \ref{t:polysg}--\ref{t:polysN} up to $i=9$.

\begin{table} 
\label{t:polysg}
\begin{tabular}{@{}rll@{}}
\hline\noalign{\smallskip}
$i$ & \multicolumn{1}{c}{$g_{i,1}$} & \multicolumn{1}{c}{$g_{i,2}$} \\
\noalign{\smallskip}\hline\noalign{\smallskip}
$1$ & $\frac{1}{4}$
\\[1ex]
$2$ & $\frac{1}{2}$
\\[1ex]
$3$ & $\frac{3}{4}+\frac{3}{8}\beta^2$ & $\frac{1}{64}$
\\[1ex]
$4$ & $1+\frac{55}{32}\beta^2$
    & $\frac{1}{16}+\frac{1}{64}\beta^2$
\\[1ex]
$5$ & $\frac{5}{4}+\frac{1251}{256}\beta^2+\frac{765}{1024}\beta^4$
    & $\frac{5}{32}+\frac{35}{256}\beta^2$
\\[1ex]
$6$ & $\frac{3}{2}+\frac{1413}{128}\beta^2+\frac{175}{32}\beta^4$
    & $\frac{5}{16}+\frac{79}{128}\beta^2+\frac{1}{16}\beta^4$
\\[1ex]
$7$ & $\frac{7}{4}+\frac{5551}{256}\beta^2+\frac{47367}{2048}\beta^4+\frac{3471}{2048}\beta^6$
    & $\frac{35}{64}+\frac{1025}{512}\beta^2+\frac{2693}{4096}\beta^4$
\\[1ex]
$8$ & $2+\frac{1239}{32}\beta^2+\frac{151301}{2048}\beta^4+\frac{140383}{8192}\beta^6$
    & $\frac{7}{8}+\frac{339}{64}\beta^2+\frac{239}{64}\beta^4+\frac{209}{1024}\beta^6$
\\[1ex]
$9$ & $\frac{9}{4}+\frac{8253}{128}\beta^2+\frac{3236383}{16384}\beta^4+\frac{3162705}{32768}\beta^6+\frac{1084959}{262144}\beta^8$
    & $\frac{21}{16}+\frac{195}{16}\beta^2+\frac{62665}{4096}\beta^4+\frac{43329}{16384}\beta^6$
\\
\noalign{\smallskip}\hline
\end{tabular}
\begin{tabular}{@{}rlll@{}}
$i$ & \multicolumn{1}{c}{$g_{i,3}$} & \multicolumn{1}{c}{$g_{i,4}$} & \multicolumn{1}{c}{$g_{i,5}$}  \\
\noalign{\smallskip}\hline\noalign{\smallskip}
$5$ & $\frac{1}{768}-\frac{1}{3072}\beta^2$
\\[1ex]
$6$ & $\frac{1}{128}\beta^2$
\\[1ex]
$7$ & $\frac{7}{256}+\frac{29}{2048}\beta^2-\frac{1}{2048}\beta^4$ & $\frac{1}{8192}-\frac{1}{16384}\beta^2$
\\[1ex]
$8$ & $\frac{7}{96}+\frac{19}{192}\beta^2+\frac{61}{8192}\beta^4$ & $\frac{1}{1024}-\frac{1}{4096}\beta^2-\frac{1}{16384}\beta^4$
\\[1ex]
$9$ & $\frac{21}{128}+\frac{1731}{4096}\beta^2+\frac{7947}{65536}\beta^4+\frac{81}{65536}\beta^6$
    & $\frac{9}{2048}+\frac{15}{16384}\beta^2-\frac{51}{65536}\beta^4$
    & $\frac{1}{81920}-\frac{3}{327680}\beta^2+\frac{1}{1310720}\beta^4$
\\
\noalign{\smallskip}\hline
\end{tabular}
\caption{Triaxiality polynomials of the transformation of $g$ in Eq.~(\protect\ref{ggp})}
\end{table}

\begin{table}
\label{t:polysl}
\begin{tabular}{@{}rll@{}}
\hline\noalign{\smallskip}
$i$ & \multicolumn{1}{c}{$\ell_{i,1}$} & \multicolumn{1}{c}{$\ell_{i,2}$} \\
\noalign{\smallskip}\hline\noalign{\smallskip}
$1$ & $\frac{1}{4}$
\\[1ex]
$2$ & $\frac{3}{4}$ & $\frac{1}{16}$
\\[1ex]
$3$ & $1+\frac{15}{32}\beta^2$ & $\frac{3}{16}$
\\[1ex]
$4$ & $\frac{5}{4}+\frac{129}{64}\beta^2$  & $\frac{25}{64}+\frac{1}{8}\beta^2$
\\[1ex]
$5$ & $\frac{3}{2}+\frac{705}{128}\beta^2+\frac{415}{512}\beta^4$ & $\frac{11}{16}+\frac{93}{128}\beta^2$
\\[1ex]
$6$ & $\frac{7}{4}+\frac{3101}{256}\beta^2+\frac{5919}{1024}\beta^4$
    & $\frac{35}{32}+\frac{649}{256}\beta^2+\frac{1179}{4096}\beta^4$
\\[1ex]
$7$ & $2+\frac{747}{32}\beta^2+\frac{24575}{1024}\beta^4+\frac{14021}{8192}\beta^6$
    & $\frac{13}{8}+\frac{439}{64}\beta^2+\frac{10191}{4096}\beta^4$
\\[1ex]
$8$ & $\frac{9}{4}+\frac{10521}{256}\beta^2+\frac{154889}{2048}\beta^4+\frac{279831}{16384}\beta^6$
    & $\frac{147}{64}+\frac{8109}{512}\beta^2+\frac{200025}{16384}\beta^4+\frac{1465}{2048}\beta^6$
\\[1ex]
$9$ & $\frac{5}{2}+\frac{4333}{64}\beta^2+\frac{1639795}{8192}\beta^4+\frac{3121345}{32768}\beta^6+\frac{262085}{65536}\beta^8$
   & $\frac{25}{8}+\frac{131}{4}\beta^2+\frac{182765}{4096}\beta^4+\frac{269247}{32768}\beta^6$
\\
\end{tabular}
\begin{tabular}{@{}rlll@{}}
\hline\noalign{\smallskip}
$i$ & \multicolumn{1}{c}{$\ell_{i,3}$} & \multicolumn{1}{c}{$\ell_{i,4}$} & \multicolumn{1}{c}{$\ell_{i,5}$} \\
\noalign{\smallskip}\hline\noalign{\smallskip}
$3$ & $\frac{1}{96}$
\\[1ex]
$4$ & $\frac{3}{64}$ & $\frac{1}{512}$
\\[1ex]
$5$ & $\frac{17}{128}+\frac{1}{32}\beta^2$ & $\frac{3}{256}$ & $\frac{1}{2560}$
\\[1ex]
$6$ & $\frac{77}{256}+\frac{117}{512}\beta^2$ & $\frac{43}{1024}+\frac{1}{128}\beta^2$ & $\frac{3}{1024}$
\\[1ex]
$7$ & $\frac{19}{32}+\frac{31}{32}\beta^2+\frac{359}{4096}\beta^4$
    & $\frac{119}{1024}+\frac{141}{2048}\beta^2$
    & $\frac{13}{1024}+\frac{1}{512}\beta^2$
\\[1ex]
$8$ & $\frac{273}{256}+\frac{6363}{2048}\beta^2+\frac{7281}{8192}\beta^4$
    & $\frac{2241}{8192}+\frac{1407}{4096}\beta^2+\frac{423}{16384}\beta^4$
    & $\frac{85}{2048}+\frac{165}{8192}\beta^2$
\\[1ex]
$9$ & $\frac{343}{192}+\frac{17095}{2048}\beta^2+\frac{82345}{16384}\beta^4+\frac{32245}{131072}\beta^6$
    & $\frac{589}{1024}+\frac{10435}{8192}\beta^2+\frac{4923}{16384}\beta^4$
    & $\frac{929}{8192}+\frac{947}{8192}\beta^2+\frac{487}{65536}\beta^4$
\\
\noalign{\smallskip}\hline
\end{tabular}
\begin{tabular}{@{}rllllllll@{}}
\noalign{\smallskip}
$i$ & \multicolumn{1}{c}{$\ell_{i,6}$} & \multicolumn{1}{c}{$\ell_{i,7}$} & \multicolumn{1}{c}{$\ell_{i,8}$} & \multicolumn{1}{c}{$\ell_{i,9}$} \\
\noalign{\smallskip}\hline\noalign{\smallskip}
$6$ & $\frac{1}{12288}$
\\[1ex]
$7$ & $\frac{3}{4096}$ & $\frac{1}{57344}$
\\[1ex]
$8$ & $\frac{61}{16384}+\frac{1}{2048}\beta^2$ & $\frac{3}{16384}$ & $\frac{1}{262144}$
\\[1ex]
$9$ & $\frac{115}{8192}+\frac{189}{32768}\beta^2$ & $\frac{35}{32768}+\frac{\beta^2}{8192}$
    & $\frac{3}{65536}$
    & $\frac{1}{1179648}$
\\
\noalign{\smallskip}\hline
\end{tabular}
\caption{Triaxiality polynomials of the transformation of $\ell$ in Eq.~(\protect\ref{llp})}
\end{table}

\begin{table}
\label{t:polysN}
\begin{tabular}{@{}rll@{}}
\hline\noalign{\smallskip}
$i$ & \multicolumn{1}{c}{$L_{i,0}$} & \multicolumn{1}{c}{$L_{i,1}$} \\
\noalign{\smallskip}\hline\noalign{\smallskip}
$1$ & $0$ & $\frac{1}{2}$
\\[1ex]
$2$ & $\frac{1}{4}$ & $\frac{1}{2}$
\\[1ex]
$3$ & $\frac{5}{8}$ & $\frac{1}{2}+\frac{3}{8}\beta^2$
\\[1ex]
$4$ & $\frac{9}{8}+\frac{27}{64}\beta^2$ & $\frac{1}{2}+\frac{21}{16}\beta^2$
\\[1ex]
$5$ & $\frac{7}{4}+\frac{35}{16}\beta^2$ & $\frac{1}{2}+\frac{387}{128}\beta^2+\frac{297}{512}\beta^4$
\\[1ex]
$6$ & $\frac{5}{2}+\frac{885}{128}\beta^2+\frac{225}{256}\beta^4$
    & $\frac{1}{2}+\frac{735}{128}\beta^2+\frac{1835}{512}\beta^4$
\\[1ex]
$7$ & $\frac{27}{8}+\frac{8775}{512}\beta^2+\frac{7155}{1024}\beta^4$
    & $\frac{1}{2}+\frac{1245}{128}\beta^2+\frac{13425}{1024}\beta^4+\frac{2307}{2048}\beta^6$
\\[1ex]
$8$ & $\frac{35}{8}+\frac{18725}{512}\beta^2+\frac{65135}{2048}\beta^4+\frac{33355}{16384}\beta^6$
    & $\frac{1}{2}+\frac{1953}{128}\beta^2+\frac{295}{8}\beta ^4+\frac{41149}{4096}\beta^6$
\\[1ex]
$9$ & $\frac{11}{2}+\frac{9009}{128}\beta^2+\frac{55583}{512}\beta^4+\frac{44825}{2048}\beta^6$
    & $\frac{1}{2}+\frac{1449}{64}\beta^2+\frac{720999}{8192}\beta^4+\frac{828675}{16384}\beta^6+\frac{321543}{131072}\beta^8$
\\
\noalign{\smallskip}\hline
\end{tabular}
\begin{tabular}{@{}rlll@{}}
$i$ & \multicolumn{1}{c}{$L_{i,2}$} & \multicolumn{1}{c}{$L_{i,3}$} & \multicolumn{1}{c}{$L_{i,4}$} \\
\noalign{\smallskip}\hline\noalign{\smallskip}
$3$ & $\frac{1}{16}$
\\[1ex]
$4$ & $\frac{3}{16}+\frac{1}{16}\beta^2$
\\[1ex]
$5$ & $\frac{3}{8}+\frac{7}{16}\beta^2$  & $\frac{1}{128}-\frac{1}{512}\beta^2$
\\[1ex]
$6$ & $\frac{5}{8}+\frac{53}{32}\beta^2+\frac{13}{64}\beta^4$ & $\frac{5}{128}+\frac{1}{512}\beta^2$
\\[1ex]
$7$ & $\frac{15}{16}+\frac{297}{64}\beta^2+\frac{473}{256}\beta^4$
    & $\frac{15}{128}+\frac{81}{1024}\beta^2-\frac{3}{2048}\beta^4$
    & $\frac{1}{1024}-\frac{1}{2048}\beta^2$
\\[1ex]
$8$ & $\frac{21}{16}+\frac{693}{64}\beta^2+\frac{9497}{1024}\beta^4+\frac{599}{1024}\beta^6$
    & $\frac{35}{128}+\frac{475}{1024}\beta^2+\frac{205}{4096}\beta^4$
    & $\frac{7}{1024}-\frac{3}{2048}\beta^2-\frac{1}{2048}\beta^4$
\\[1ex]
$9$ & $\frac{7}{4}+\frac{357}{16}\beta^2+\frac{69863}{2048}\beta^4+\frac{3487}{512}\beta^6$
    & $\frac{35}{64}+\frac{3575}{2048}\beta^2+\frac{20229}{32768}\beta^4+\frac{45}{4096}\beta^6$
    & $\frac{7}{256}+\frac{35}{4096}\beta^2-\frac{11}{2048}\beta^4$
\\
\noalign{\smallskip}\hline
\noalign{\smallskip}
\end{tabular}
\begin{tabular}{@{}rl@{}}
\noalign{\smallskip}
 & $L_{9,5}=\frac{1}{8192}-\frac{3}{32768}\beta^2+\frac{1}{131072}\beta^4$
\\
\noalign{\smallskip}\hline
\end{tabular}
\caption{Triaxiality polynomials of the transformation of $N$ in Eq.~(\protect\ref{NNp})}
\end{table}

It should be said that the Lie series solution has been computed up to an order much higher than that required in common astronomical problems. But the polynomials provided in Tables \ref{t:polyham}--\ref{t:polysN} can be useful for those interested in checking their own implementations of the method. Besides, results of this paper are not limited to astronomical studies and also apply to the general problem of rigid body rotation close to the axis of maximum inertia, where such orders may be required. The computation of even higher orders, if needed, is easily achieved in the framework provided by Deprit's formulation of  the Lie transforms algorithm \cite{Deprit1969}.

\section{Comparison with Kinoshita's expansions}

As a way of easing development of the disturbing function in the action angle variables of the torque free motion, Kinoshita provided different expansions that apply for short- and long-axis modes \cite{Kinoshita1992}. These expansions are customarily used in the study of the rotation of solar system bodies (see, for instance, \cite{Kinoshita1977,GetinoFerrandiz1991,SouchayKinoshitaNakaiRoux2003,CottereauSouchayAljbaae2010}).
\par

For the present case of concern of rotations in the short-axis-mode, it is found the relation
\begin{equation}\label{jotita}
1-\cos{j}=1-\frac{N}{M}=\frac{1-\beta}{\sqrt{1-\beta^2}}\,\frac{L}{G}=(1-\beta)\,\delta,
\end{equation}
where $j$ corresponds to the minimum value of $J$, which is one of the integration constants of Kinoshita's solution\cite{Kinoshita1992}. Then, replacing Eq.~(\ref{NNp}) into Eq.~(\ref{jotita}) and working to the order of $j^2$ results in
\begin{equation}\label{dj}
\frac{j^2}{2}=
(1-\beta)\,\delta'=\sqrt{\frac{1-\beta}{1+\beta}}\,\frac{L'}{G'}.
\end{equation}
\par

Equation (\ref{dj}) allows to check that Kinoshita's expansions for short-axis-mode \cite{Kinoshita1992} are in agreement with the lower order terms of the Lie transforms theory of this paper. Thus, for instance, the secular frequency of $\ell$ is obtained from the corresponding Hamilton equation of the averaged Hamiltonian (\ref{averaged}) as
\[
n_{\ell'}=\frac{\partial\mathcal{T}}{\partial{L}}.
\]
Up to the second order, which means making $i=1$ in Eq.~(\ref{averaged}), it is found
\begin{eqnarray*}
n_{\ell'} &=& \frac{G'}{C}\,\alpha\left(\sqrt{1-\beta ^2}-\frac{L'}{G'}-\frac{3}{4}\,\frac{\beta^2}{\sqrt{1-\beta^2}}\,\frac{L'^2}{G'^2}\right).
\end{eqnarray*}
Then, using Eq.~(\ref{dj}),
\begin{eqnarray*}
n_{\ell'} &=& \frac{G'}{C}\,\alpha\sqrt{1-\beta ^2}\left[1-\frac{1}{2(1-\beta)}\,j^2-\frac{3\beta^2}{16 (1-\beta)^2}\,j^4\right]
\end{eqnarray*}
which, up to the order of $j^2$, is in agreement with the expression for $n_{\tilde{l}}$ in Eq.~(39) of \cite{Kinoshita1992}.
\par

Proceeding analogously, the secular frequency of $g$ is found to be
\[
n_{g'}=\frac{G'}{C}\left[1+\frac{\alpha}{2}(1+\beta)\,j^2+\frac{\alpha}{32}\,(1+\beta)\frac{\beta^2}{(1-\beta)^2}\,j^6\right].
\]
Then, taking into account Eq.~(\ref{mug}) and up to the order of $j^2$, 
\[
n_{g'}+n_{\ell'}=\frac{G'}{C}+\frac{G'}{C}\,\alpha\left[
1-\left(1-\sqrt{1-\beta ^2}\right) \left(1-\frac{1}{2}\sqrt{\frac{1+\beta}{1-\beta}}\,j^2\right)
\right],
\]
which is in agreement with the expansion for $n_{\tilde{g}}$ in Eq.~(40) of \cite{Kinoshita1992}.
\par

In what respects to the variables, the procedure for checking Kinoshita's expansions is slightly more involved because the argument $\ell$ of the circular functions in the initial transformation to action-angle variables, Eqs.~(\ref{scnu})--(\ref{NL}), is given by the Lie series in Eq.~(\ref{llp}), and therefore must be expanded in power series of $j$. Thus, for instance, from Eqs.~(\ref{NL}), (\ref{llp}), (\ref{NNp}), Tables \ref{t:polysl} and \ref{t:polysN}, and Eq.~(\ref{dj})
\begin{eqnarray*}
N &=& G'-L'\left[1+\frac{\beta}{1-\beta}\,\frac{j^2}{4}\cos2\ell'+\mathcal{O}(j^4)\right] \\
&& \times\frac{1}{\sqrt{1-\beta^2}}\left\{1+\beta\cos\left[2\ell'-\frac{\beta}{1-\beta}\,\frac{j^2}{4}\sin2\ell'+\mathcal{O}(j^4)\right]\right\},
\end{eqnarray*}
which after expansion in powers of $j$ results in
\[
N=M\left(1-\frac{1}{2}j^2-\frac{\beta}{1-\beta}\,j^2\cos^2\ell'\right)+\mathcal{O}(j^4),
\]
in agreement with Eq.~(36) of \cite{Kinoshita1992}.

\section{Long-axis-mode}

As it is well-known, SAM and LAM rotations are formulated analogously by the simple procedure of interchanging the moments of maximum and minimum inertia, and properly labeling the body axes. Because of this interchange of $A$ and $C$, Eq.~(\ref{alfabeta}) provides  new inertia parameters given by
\begin{equation} \label{alfabetastar}
\alpha^*\,(1+\beta^*)=\frac{A}{C}-1,
\qquad
\alpha^*\,(1-\beta^*)=\frac{A}{B}-1.
\end{equation}
A thorough discussion on inertia parameters for the free rigid body problem can be found in \cite{HitzlBreakwell1971}, where it is pointed out the relation
\begin{equation} \label{betastar}
\beta^*=\frac{1-\beta}{1+3\beta},
\end{equation}
which is later used by Kinoshita in his derivation of LAM analytical expansions \cite{Kinoshita1992}.
\par

Thus, Eq.~(\ref{hacheAsin}) is reformulated as
\begin{equation}\label{fuku}
\mathcal{H}_0=\frac{M^2}{2A}\left[1+\alpha^*\sin^2J^*\left(1+\beta^*\cos2\nu^*\right)\right],
\end{equation}
where $J^*$ is the inclination angle between the plane $\Pi$ perpendicular to the body's $x$ axis and the invariable plane, and $\nu^*$ is the angle between the ascending node of $\Pi$ in the invariable plane and the body's $y$ axis. 
Therefore, the SAM discussion above applies also to LAM except for a change in the sign of the inertia parameter $\beta^*$. 

\section{Conclusions}

An alternative solution to the torque free motion of a rigid body has been constructed, which applies to the particular case of rotation close to the principal axis of maximum inertia and is valid for bodies with any triaxiality. For this specific case, the torque free motion Hamiltonian is rearranged as a perturbation problem, the main part of which is formally the same as the full problem with the proviso of using half the inclination angle between the equatorial plane of the rigid body and the plane perpendicular to the angular momentum vector, instead of the full angle. Therefore, in those cases in which this inclination angle is small, the main part of the short-axis-mode rotation may provide by itself an accurate description of the torque free motion by means of a solution that only involves trigonometric functions. Besides, after reformulating the main part of the Hamiltonian in action-angle variables, the new solution fits naturally in a general perturbation scheme, thus paving the way in the computation of perturbation solutions to rigid body rotation under external torques. Results in this direction are in progress and, hopefully, will be published elsewhere.

\section{Acknowledgment}
Support is acknowledged from projects AYA 2009-11896 and AYA 2010-18796 of the Government of Spain.


\begin{thebibliography}{99}

\bibitem{Kinoshita1977}
H.~{Kinoshita}, ``{Theory of the rotation of the rigid earth},''  {\em
  Celestial Mechanics}, Vol.~15, No.~3, 1977, pp.~277--326, 10.1007/BF01228425.

\bibitem{GetinoFerrandiz1991}
J.~{Getino} and J.~M. {Ferr{\'a}ndiz}, ``{A Hamiltonian Theory for an Elastic
  Earth - First Order Analytical Integration},''  {\em Celestial Mechanics and
  Dynamical Astronomy}, Vol.~51, No.~1, 1991, pp.~35--65, 10.1007/BF02426669.

\bibitem{CicaloScheeres2010}
S.~{Cical{\`o}} and D.~J. {Scheeres}, ``{Averaged rotational dynamics of an
  asteroid in tumbling rotation under the YORP torque},''  {\em Celestial
  Mechanics and Dynamical Astronomy}, Vol.~106, No.~4, 2010, pp.~301--337,
  10.1007/s10569-009-9249-7.

\bibitem{Sadov1970}
Y.~A. Sadov, ``The Action-Angles Variables in the Euler-Poinsot Problem,''
  {\em PMM-Journal of Applied Mathematics and Mechanics}, Vol.~34, No.~5, 1970,
  pp.~922--925.

\bibitem{Kinoshita1972}
H.~{Kinoshita}, ``{First-Order Perturbations of the Two Finite Body Problem},''
   {\em Publications of the Astronomical Society of Japan}, Vol.~24, 1972,
  pp.~423--457.

\bibitem{Kozlov1974}
V.~V. {Kozlov}, ``{La G\'eom\'etrie des variables action-angle dans le
  probl\`eme d'Euler-Poinsot},''  {\em Vestnik Moskovskogo Universiteta. Seriya
  I. Matematika, Mekhanika}, Vol.~5, 1974, pp.~74--79.
\newblock (in Russian).

\bibitem{LaraFerrer2013}
M.~{Lara} and S.~{Ferrer}, ``{Closed Form Perturbation Solution of a Fast
  Rotating Triaxial Satellite under Gravity-Gradient Torque},''  {\em Cosmic
  Research}, Vol.~51, No.~4, 2013, pp.~289--303.

\bibitem{Arnold1989}
V.~I. Arnold, {\em Mathematical Methods of Classical Mechanics}.
\newblock New York: Springer-Verlag, 2nd~ed., 1989.

\bibitem{LaraFukushimaFerrer2010}
M.~{Lara}, T.~{Fukushima}, and S.~{Ferrer}, ``{First-order rotation solution of
  an oblate rigid body under the torque of a perturber in circular orbit},''
  {\em Astronomy \& Astrophysics}, Vol.~519, 2010, p.~A1,
  10.1051/0004-6361/200913880.

\bibitem{Barkin1998}
Y.~V. {Barkin}, ``{Unperturbed Chandler Motion and Perturbation Theory of the
  Rotation Motion of Deformable Celestial Bodies},''  {\em Astronomical and
  Astrophysical Transactions}, Vol.~17, No.~3, 1998, pp.~179--219,
  10.1080/10556799808232092.

\bibitem{Chernousko1963}
F.~L. Chernous'ko, ``On the motion of a satellite about its center of mass
  under the action of gravitational moments,''  {\em PMM-Journal of Applied
  Mathematics and Mechanics}, Vol.~27, No.~3, 1963, pp.~708--722.

\bibitem{HitzlBreakwell1971}
D.~L. {Hitzl} and J.~V. {Breakwell}, ``{Resonant and non-resonant
  gravity-gradient perturbations of a tumbling tri-axial satellite.},''  {\em
  Celestial Mechanics}, Vol.~3, No.~5, 1971, pp.~346--383, 10.1007/BF01231806.

\bibitem{Kinoshita1992}
H.~{Kinoshita}, ``{Analytical expansions of torque-free motions for short and
  long axis modes},''  {\em Celestial Mechanics and Dynamical Astronomy},
  Vol.~53, No.~4, 1992, pp.~365--375, 10.1007/BF00051817.

\bibitem{SouchayFolgueiraBouquillon2003}
J.~{Souchay}, M.~{Folgueira}, and S.~{Bouquillon}, ``{Effects of the
  Triaxiality on the Rotation of Celestial Bodies: Application to the Earth,
  Mars and Eros},''  {\em Earth Moon and Planets}, Vol.~93, No.~2, 2003,
  pp.~107--144, 10.1023/B:MOON.0000034505.79534.01.

\bibitem{GetinoEscapaMiguel2010}
J.~{Getino}, A.~{Escapa}, and D.~{Miguel}, ``{General Theory of the Rotation of
  the Non-rigid Earth at the Second Order. I. The Rigid Model in Andoyer
  Variables},''  {\em Astronomical Journal}, Vol.~139, No.~5, 2010,
  pp.~1916--1934, 10.1088/0004-6256/139/5/1916.

\bibitem{CottereauSouchayAljbaae2010}
L.~{Cottereau}, J.~{Souchay}, and S.~{Aljbaae}, ``{Accurate free and forced
  rotational motions of rigid Venus},''  {\em Astronomy \& Astrophysics},
  Vol.~515, 2010, p.~A9, 10.1051/0004-6361/200913785.

\bibitem{SouchayKinoshitaNakaiRoux2003}
J.~{Souchay}, H.~{Kinoshita}, H.~{Nakai}, and S.~{Roux}, ``{A precise modeling
  of Eros 433 rotation},''  {\em Icarus}, Vol.~166, No.~2, 2003, pp.~285--296,
  10.1016/j.icarus.2003.08.018.

\bibitem{SouchayBouquillon2005}
J.~{Souchay} and S.~{Bouquillon}, ``{The high frequency variations in the
  rotation of Eros},''  {\em Astronomy \& Astrophysics}, Vol.~433, No.~1, 2005,
  pp.~375--383, 10.1051/0004-6361:20035780.

\bibitem{Zanardi1986}
M.~C. {Zanardi}, ``{Study of the terms of coupling between rotational and
  translational motions},''  {\em Celestial Mechanics}, Vol.~39, No.~1, 1986,
  pp.~147--158, 10.1007/BF01230847.

\bibitem{FerrandizSansaturio1989}
J.-M. {Ferrandiz} and M.-E. {Sansaturio}, ``Elimination of the nodes when the
  satellite is a non spherical rigid body,''  {\em Celestial Mechanics and
  Dynamical Astronomy}, Vol.~46, No.~4, 1989, pp.~307--320, 10.1007/BF00051485.

\bibitem{FerrerLara2010}
S.~{Ferrer} and M.~{Lara}, ``{Integration of the Rotation of an Earth-like Body
  as a Perturbed Spherical Rotor},''  {\em The Astronomical Journal}, Vol.~139,
  No.~5, 2010, pp.~1899--1908, 10.1088/0004-6256/139/5/1899.

\bibitem{LaraFukushimaFerrer2011}
M.~{Lara}, T.~{Fukushima}, and S.~{Ferrer}, ``{Ceres' rotation solution under
  the gravitational torque of the Sun},''  {\em Monthly Notices of the Royal
  Astronomical Society}, Vol.~415, No.~1, 2011, pp.~461--469,
  10.1111/j.1365-2966.2011.18717.x.

\bibitem{Fukushima2008}
T.~{Fukushima}, ``{Simple, Regular, and Efficient Numerical Integration of
  Rotational Motion},''  {\em Astronomical Journal}, Vol.~135, No.~6, 2008,
  pp.~2298--2322, 10.1088/0004-6256/135/6/2298.

\bibitem{Hori1966}
G.~{Hori}, ``{Theory of General Perturbation with Unspecified Canonical
  Variable},''  {\em Publications of the Astronomical Society of Japan},
  Vol.~18, No.~4, 1966, p.~287.

\bibitem{Deprit1969}
A.~{Deprit}, ``Canonical transformations depending on a small parameter,''
  {\em Celestial Mechanics}, Vol.~1, No.~1, 1969, pp.~12--30,
  10.1007/BF01230629.

\bibitem{Whittaker1917}
E.~T. {Whittaker}, {\em {A Treatise on the Analytical Dynamics of Particles and
  Rigid Bodies}}.
\newblock Cambridge University Press, 2nd~ed., Feb. 1917.

\bibitem{Golubev1960}
V.~Golubev, {\em Lectures on Integration of the Equations of Motion of a Rigid
  Body about a Fixed Point}.
\newblock Jerusalem: Israel Program for Scientific Translations, S.~Monson,
  1960.

\bibitem{Andoyer1923}
M.~H. Andoyer, {\em Cours de M\'ecanique C\'eleste}.
\newblock Paris: Gauthier-Villars et cie, 1923.

\bibitem{DepritAJP1967}
A.~{Deprit}, ``Free Rotation of a Rigid Body Studied in the Phase Space,''
  {\em American Journal of Physics}, Vol.~35, 1967, pp.~424--428.

\bibitem{FerrerLara2010b}
S.~{Ferrer} and M.~{Lara}, ``{Families of Canonical Transformations by
  Hamilton-Jacobi-Poincar\'e Equation. Application to Rotational and Orbital
  Motion},''  {\em Journal of Geometric Mechanics}, Vol.~2, No.~3, 2010,
  pp.~223--241, 10.3934/jgm.2010.2.223.

\bibitem{DehantdeVironKaratekinVanHoolst2006}
V.~{Dehant}, O.~{de Viron}, O.~{Karatekin}, and T.~{van Hoolst}, ``{Excitation
  of Mars polar motion},''  {\em Astronomy \& Astrophysics}, Vol.~446, No.~1,
  2006, pp.~345--355, 10.1051/0004-6361:20053825.

\bibitem{Fukushima2006}
T.~Fukushima, ``Efficient integration of torque-free rotation by energy scaling
  method,''  {\em Proceedings of the Journ\'ees Syst\`emes de R\'ef\'erence
  Spatio-Temporels 2005} (A.~{Brzezinski}, N.~{Capitaine}, and B.~{Kolaczek},
  eds.), Space Research Centre PAS, Warsaw, Poland, 2006.

\bibitem{NewhallWilliams1997}
X.~X. {Newhall} and J.~G. {Williams}, ``{Estimation of the Lunar Physical
  Librations},''  {\em Celestial Mechanics and Dynamical Astronomy}, Vol.~66,
  No.~1, 1996, pp.~21--30, 10.1007/BF00048820.

\end{thebibliography}
\end{document}